\documentclass{article}

\usepackage{amssymb,amsmath,color}
\usepackage{amsthm}
\usepackage{url}
\usepackage{tikz}

\definecolor{red}{rgb}{1,0,0}
\definecolor{blue}{rgb}{.2,.2,.8}

\newtheorem{theorem}{Theorem}[section]
\newtheorem{corollary}[theorem]{Corollary}

\newtheorem{conjecture}{Conjecture}

\theoremstyle{definition}

\begin{document}

\title{$q$-Series congruences involving statistical mechanics partition functions in regime III and IV of Baxter's solution of the hard-hexagon model}
\author{Mircea Merca
	\\ 
	\footnotesize Department of Mathematics, University of Craiova, Romania\\
	\footnotesize Academy of Romanian Scientists, Bucharest, Romania\\
	\footnotesize mircea.merca@profinfo.edu.ro
}
\date{}
\maketitle

\begin{abstract} 
	For each $s\in\{2,4\}$, the generating function of $R_s(n)$, the number of partitions of $n$ into odd parts or congruent to $0$, $\pm s\pmod {10}$, arises naturally in regime III of Rodney Baxter's solution of the hard-hexagon model of statistical mechanics.  
	For each $s\in\{1,3\}$, the generating function of $R^*_s(n)$, the number of partitions of $n$ into parts not congruent to $0$, $\pm s\pmod {10}$ and $10-2s \pmod {20}$, arises naturally in regime IV of Rodney Baxter's solution of the hard-hexagon model of statistical mechanics. 
	In this paper, we  investigate the parity of $R_s(n)$ and $R^*_s(n)$, providing  new parity results involving sums of partition numbers $p(n)$ and squares in arithmetic progressions.
\end{abstract}

\section{Introduction}

A partition of a positive integer $n$ is a sequence of positive integers whose sum is $n$. The order of the summands is unimportant when writing the partitions of $n$, but for consistency, a partition of $n$ will be written with the summands in a nonincreasing order \cite{Andrews98}. As usual, we denote by $p(n)$ the number of the partitions of $n$. For example, we have $p(5)=7$ because the partitions of $5$ are given as:
$$ 5,\ 4+1,\ 3+2,\ 3+1+1,\ 2+2+1,\ 2+1+1+1,\ 1+1+1+1+1.$$
The fastest algorithms for enumerating all the partitions of an integer have recently been presented by Merca \cite{Merca12,Merca13}.

The partition function $p(n)$ may be defined by the generating function
\begin{equation*}\label{GFP}
\sum_{n=0}^{\infty}p(n)q^n=\frac{1}{(q;q)_{\infty}}.
\end{equation*}
Here and throughout this paper, we use the following customary $q$-series notation:
$$(a;q)_0=1,\qquad (a;q)_n=\prod_{k=0}^{n-1}(1-aq^k),\qquad (a;q)_\infty=\prod_{k=0}^{\infty}(1-aq^k).$$
Because the infinite product $(a;q)_{\infty}$ diverges when $a\neq 0$ and $|q|\geqslant 1$, whenever $(a;q)_{\infty}$ appears in a formula, we shall assume $|q|<1$
and we shall use the compact notation
$$(a_1,a_2\ldots,a_n;q)_{\infty}=(a_1;q)_{\infty}(a_2;q)_{\infty}\cdots(a_n;q)_{\infty}.$$

For $s\in\{2,4\}$, we denote by $R_s(n)$ the number of partitions of $n$ into odd parts or congruent to $0$, $\pm s\pmod {10}$. Elementary techniques in the theory of partitions give the following generating function for $R_s(n)$:
$$\sum_{n=0}^{\infty} R_s(n) q^n = \frac{1}{(q;q^2)_{\infty}(q^s,q^{10-s};q^{10})_{\infty}}.$$
The generating function for $R_s(n)$ arises naturally in regime III of Rodney Baxter's solution of the hard-hexagon model of  statistical mechanics, appearing in the following $q$-identities of Rogers \cite{Andrews81,Andrews07}:
$$\sum_{n=0}^{\infty} \frac{(-q;q)_nq^{n(3n+1)/2}}{(q;q)_{2n+1}} = \frac{G(q^2)}{(q;q^2)_{\infty}},$$
and
$$\sum_{n=0}^{\infty} \frac{(-q;q)_nq^{3n(n+1)/2}}{(q;q)_{2n+1}} = \frac{H(q^2)}{(q;q^2)_{\infty}},$$
where
$$G(q):=\sum_{n=0}^{\infty}\frac{q^{n^2}}{(q;q)_n}=\frac{1}{(q,q^4;q^5)_{\infty}}$$
and
$$H(q):=\sum_{n=0}^{\infty}\frac{q^{n^2+n}}{(q;q)_n}=\frac{1}{(q^2,q^3;q^5)_{\infty}}$$
are the Rogers-Ramanujan functions. %

For $s\in\{1,3\}$, we denote by $R^*_s(n)$ the number of partitions of $n$ into parts not congruent to $0$, $\pm s\pmod {10}$ and $10-2s \pmod {20}$. Elementary techniques in the theory of partitions give the following generating function for $R^*_s(n)$:
$$\sum_{n=0}^{\infty} R^*_s(n) q^n = \frac{(q^s,q^{10-s},q^{10};q^{10})_\infty\,(q^{10-2s},q^{10+2s};q^{20})_\infty}{(q;q)_{\infty}}.$$
The generating function for $R^*_s(n)$ arises naturally in regime IV of Rodney Baxter's solution of the hard-hexagon model of  statistical mechanics, appearing in the following $q$-identities of Rogers \cite{Andrews81,Andrews07}:
$$\sum_{n=0}^{\infty} \frac{q^{n(n+1)}}{(q;q)_{2n}} = \frac{(q,q^{9},q^{10};q^{10})_\infty\,(q^{8},q^{12};q^{20})_\infty}{(q;q)_{\infty}}$$
and
$$\sum_{n=0}^{\infty} \frac{q^{n(n+1)}}{(q;q)_{2n+1}} = \frac{(q^3,q^{7},q^{10};q^{10})_\infty\,(q^{4},q^{16};q^{20})_\infty}{(q;q)_{\infty}}.$$

In this paper, we shall provide the following $q$-series congruences.

\begin{theorem}\label{T1}
	Let $s$ be a positive integer.
	\begin{enumerate}
		\item For $s\in\{2,4\}$,
		$$
		\sum_{n=0}^{\infty} \frac{(-q;q)_n q^{n(3n+s-1)/2}}{(q;q)_{2n+1}}
		\equiv \sum_{120n+(3s-5)^2\text{ square}} q^n \pmod 2.
		$$	
		\item For $s\in\{1,3\}$,
		$$
		\sum_{n=0}^{\infty} \frac{q^{n(n+1)}}{(q;q)_{2n+(s-1)/2}}
		\equiv \sum_{40n+s^2\text{ square}} q^n \pmod 2.
		$$
	\end{enumerate}
\end{theorem}

As a consequence of Theorem \ref{T1}, we immediately deduce the following parity result involving sums of partition numbers $p(n)$ and squares in arithmetic progressions.

\begin{corollary}\label{C2}
	Let $n$ be a nonnegative integer. 
	\begin{enumerate}
		\item For $s\in\{2,4\}$,
		$$
		\sum_{20k+(s-1)^2\text{ square}} p(n-k) \equiv 1 \pmod 2
		$$
		if and only if $120n+(3s-5)^2$ is a square.
		\item 	For $s\in\{1,3\}$,
		$$
		\sum_{15k+(s+1)^2/4\text{ square}} p(n-k) \equiv 1 \pmod 2
		$$
		if and only if $40n+s^2$ is a square.
	\end{enumerate}
\end{corollary}

 Questions regarding the parity of sums of partition numbers for square values in given arithmetic progressions have been studied recently \cite{BM}. 
 As we can see in \cite{BM}, the cases $s=1$ and $s=2$ of Corollary \ref{C2} are known.
 The following conjecture is also known: the statement 
\begin{equation*}\label{main}
\sum_{ak+1\text{ square}} p(n-k) \equiv 1 \pmod 2 \quad\mbox{if and only if}\quad bn+1 \mbox{ is a square}.
\end{equation*}
is true if and only if $$(a,b) \in S:=\{(6,8), (8,12), (12, 24), (15, 40), (16, 48), (20, 120), (21, 168)\}.$$
This paper opens new possibilities for research in this area.

The organization of the paper is as follows. We will first prove Theorem \ref{T1}
in Section \ref{S2} considering a truncated theta identity of Gauss \cite{Andrews18}. 
In Section \ref{S3}, we will provide a proof of Corollary \ref{C2}
considering a decomposition of $R_s(n)$ in terms of Euler partition function $p(n)$.
In Section \ref{S4}, we will introduce new open problems involving the partition function $R_s(n)$, $s\in\{2,4\}$ and $R^*_s(n)$, $s\in\{1,3\}$.  
  
\section{Proof of Theorem \ref{T1}}
\label{S2}

Watson's quintuple product identity \cite{Carlitz,Subbarao} states that 
\begin{equation}\label{WQPI}
\sum_{n=-\infty}^{\infty}z^{3n}q^{n(3n-1)/2}(1-zq^n)
=(q,z,q/z;q)_{\infty}(qz^2,q/z^2;q^2)_{\infty}.
\end{equation}
By this identity, with $q$ replaced by $-q^{5}$ and $z$ replaced by $q^2$, we get
\allowdisplaybreaks{
	\begin{align*}
	&\sum_{n=-\infty}^{\infty}(-1)^{n(3n-1)/2}\left( q^{n(15n+7)/2}+q^{(3n-2)(5n-1)/2}\right) \\
	&\qquad\qquad=(q^2,-q^3,-q^5;-q^5)_{\infty}\cdot (-q,-q^9;q^{10})_{\infty}\\
	&\qquad\qquad=(q^2,-q^3,-q^5;-q^5)_{\infty}\cdot \frac{(q^2,q^{18};q^{20})_{\infty}}{(q,q^9;q^{10})_{\infty}}\\
	&\qquad\qquad=(-q^3,-q^5,-q^7;q^{10})_{\infty}\cdot (q^2,q^8,q^{10};q^{10})_{\infty}\cdot \frac{(q^2,q^{18};q^{20})_{\infty}}{(q,q^9;q^{10})_{\infty}}\\
	&\qquad\qquad=\frac{(q^6,q^{10},q^{14};q^{20})_{\infty}}{(q^3,q^5,q^7;q^{10})_{\infty}}\cdot \frac{(q^2,q^4,q^6,q^8,q^{10};q^{10})_{\infty}}{(q^4,q^6;q^{10})_{\infty}}\cdot \frac{(q^2,q^{18};q^{20})_{\infty}}{(q,q^9;q^{10})_{\infty}}\\
	&\qquad\qquad=\frac{(q^2,q^6,q^{10},q^{14},q^{18};q^{20})_{\infty}}{(q,q^3,q^5,q^7,q^9;q^{10})_{\infty}}\cdot \frac{(q^2;q^{2})_{\infty}}{(q^4,q^6;q^{10})_{\infty}}\\
	&\qquad\qquad=\frac{(q^2;q^{4})_{\infty}}{(q;q^{2})_{\infty}}\cdot \frac{(q^2;q^{2})_{\infty}}{(q^4,q^6;q^{10})_{\infty}}\\
	&\qquad\qquad=\frac{1}{(q;q^{2})_{\infty}(q^4,q^6;q^{10})_{\infty}}\cdot \frac{(q^2;q^2)_\infty}{(-q^2;q^2)_\infty}.
	\end{align*}}
In a similar way, letting $q\to\-q^5$ and setting $z=-q$ in \eqref{WQPI}, we obtain
\begin{align*}
& \sum_{n=-\infty}^{\infty}(-1)^{n(n-1)/2}\left( q^{n(15n+1)/2}+q^{(3n-1)(5n-2)/2}\right) \\
& \qquad\qquad = \frac{1}{(q;q^{2})_{\infty}(q^2,q^8;q^{10})_{\infty}}\cdot \frac{(q^2;q^2)_\infty}{(-q^2;q^2)_\infty}.
\end{align*}
Thus, for $s\in\{2,4\}$, we deduce that
\begin{align}
& \sum_{n=-\infty}^{\infty}(-1)^{n((s-1)n-1)/2}\left( q^{n(15n+3s-5)/2}+q^{(3n-s/2)(5n-3+s/2)/2}\right) \nonumber \\
& \qquad\qquad= \frac{1}{(q;q^{2})_{\infty}(q^s,q^{10-s};q^{10})_{\infty}}\cdot \frac{(q^2;q^2)_\infty}{(-q^2;q^2)_\infty}. \label{eq4.1}
\end{align}

On the other hand, the Jacobi triple product identity says that
(cf.\ \cite[Eq.~(1.6.1)]{GaRaAA})
\begin{equation} \label{eq:JTP} 
\sum_{n=-\infty}^\infty (-z)^n q^{n(n-1)/2}
=
(z,q/z,q;q)_\infty.
\end{equation}
By this identity, with $q$ replaced by $-q^{5}$ and $z$ replaced by $-q$, we derive
\allowdisplaybreaks{
	\begin{align*}
	& \sum_{n=-\infty}^{\infty} (-1)^{n(n-1)/2}q^{n(5n-3)/2} \\ 
	& \qquad\qquad =(-q,q^4,-q^5;-q^5)_{\infty} \\
	& \qquad\qquad =(q^4,q^6,q^{10};q^{10})_{\infty}(-q,-q^5,-q^{9};q^{10})_{\infty} \\
	& \qquad\qquad =(q^4,q^6,q^{10};q^{10})_{\infty} \frac{(q^2,q^{10},q^{18};q^{20})_{\infty}}{ (q,q^5,q^{9};q^{10})_{\infty} } \\
	& \qquad\qquad =\frac { (q^2,q^4,q^6,q^8,q^{10};q^{10})_{\infty} } {(q,q^2,q^5,q^8,q^9;q^{10})_{\infty}} \cdot \frac{(q^2,q^6,q^{10},q^{14},q^{18};q^{20})_{\infty}}{ (q^6,q^{14};q^{20})_{\infty} } \\
	& \qquad\qquad =\frac { (q^2;q^{2})_{\infty} } {(q,q^2,q^5,q^8,q^9;q^{10})_{\infty}} \cdot \frac{(q^2;q^{4})_{\infty}}{ (q^6,q^{14};q^{20})_{\infty} } \\
	& \qquad\qquad =\frac { (q^3,q^4,q^6,q^7,q^{10};q^{10})_{\infty} } {(q;q)_{\infty}} \cdot \frac{(q^2,q^2,q^4;q^{4})_{\infty}}{ (q^6,q^{14};q^{20})_{\infty} } \\	
	& \qquad\qquad =\frac { (q^3,q^4,q^7,q^{10},q^{13},q^{16},q^{17},q^{20};q^{20})_{\infty} } {(q;q)_{\infty}} \cdot (q^2,q^2,q^4;q^{4})_{\infty} \\
	& \qquad\qquad =\frac { (q^3,q^7,q^{10};q^{10})_{\infty}\,(q^4,q^{16};q^{20})_{\infty} } {(q;q)_{\infty}} \cdot \frac{(q^2;q^{2})_{\infty}}{(-q^2;q^{2})_{\infty}}.
	\end{align*}
}
In a similar way, letting $q\to -q^5$ and setting $z=-q^3$ in \eqref{eq:JTP}, we obtain
$$\sum_{n=-\infty}^{\infty} (-1)^{n(n-1)/2}q^{n(5n+1)/2}
= \frac { (q,q^9,q^{10};q^{10})_{\infty}\,(q^8,q^{12};q^{20})_{\infty} } {(q;q)_{\infty}} \cdot \frac{(q^2;q^{2})_{\infty}}{(-q^2;q^{2})_{\infty}}.$$
Thus, for $s\in\{1,3\}$, we deduce that
\begin{align}
& \sum_{n=-\infty}^{\infty} (-1)^{n(n+s)/2}q^{n(5n-s)/2} \nonumber \\
& \qquad\qquad = \frac { (q^2,q^{10-s},q^{10};q^{10})_{\infty}\,(q^{10-2s},q^{10+2s};q^{20})_{\infty} } {(q;q)_{\infty}} \cdot \frac{(q^2;q^{2})_{\infty}}{(-q^2;q^{2})_{\infty}}. \label{eq4.2}
\end{align}

The following theta identity is often attributed to Gauss \cite[p. 23, eqs. (2.2.12)]{Andrews98}:
\begin{equation*}\label{eq:Gauss}
1+2\sum_{n=1}^{\infty} (-1)^n q^{n^2} = \frac {(q;q)_\infty} {(-q;q)_\infty}.
\end{equation*}
In \cite{Andrews18}, the authors considered this theta identity and proved the following truncated form:
\begin{align*}
& \frac{(-q;q)_{\infty}} {(q;q)_{\infty}} \left(1 + 2 \sum_{j=1}^{k} (-1)^j q^{j^2} \right) \\
& \qquad\qquad = 1+2 (-1)^k \frac{(-q;q)_k}{(q;q)_k} \sum_{j=0}^{\infty} \frac{q^{(k+1)(k+j+1)}(-q^{k+j+2};q)_{\infty}}{(1-q^{k+j+1})(q^{k+j+2};q)_{\infty}}.
\end{align*}
By this identity, with $q$ replaced by $q^2$, we obtain
\begin{align}
& (-1)^k \left( \frac{(-q^2;q^2)_{\infty}} {(q^2;q^2)_{\infty}} \left(1 + 2 \sum_{j=1}^{k} (-1)^j q^{2j^2} \right) -1\right) \nonumber \\
& \qquad\qquad = 2\frac{(-q^2;q^2)_k}{(q^2;q^2)_k} \sum_{j=0}^{\infty} \frac{q^{2(k+1)(k+j+1)}(-q^{2(k+j+2)};q^2)_{\infty}}{(1-q^{2(k+j+1)})(q^{2(k+j+2)};q^2)_{\infty}}. \label{eq4.3}
\end{align}
Multiplying both sides of this identity by \eqref{eq4.1}, we obtain
\begin{align}
& \bigg(1 + 2 \sum_{n=1}^{k} (-1)^n q^{2n^2} \bigg) \sum_{n=0}^{\infty} \frac{(-q;q)_n q^{n(3n+s-1)/2}}{(q;q)_{2n+1}} \nonumber \\ 
& \qquad -\sum_{n=-\infty}^{\infty}(-1)^{n((s-1)n-1)/2}\left( q^{n(15n+3s-5)/2}+q^{(3n-s/2)(5n-3+s/2)/2}\right) \nonumber \\
& =2(-1)^k \frac{(-q^2;q^2)_k}{(q^2;q^2)_k} \sum_{n=k+1}^{\infty} \frac{q^{2n(k+1)}(-q^{2(n+1)};q^2)_{\infty}}{(1-q^{2n})(q^{2(n+1)};q^2)_{\infty}}\times \nonumber \\
& \qquad \times \sum_{n=-\infty}^{\infty}(-1)^{n((s-1)n-1)/2}\left( q^{n(15n+3s-5)/2}+q^{(3n-s/2)(5n-3+s/2)/2}\right).\label{Id1}
\end{align}
For each $s\in\{2,4\}$, it is elementary to see that squares that are congruent to $(3s-5)^2$ modulo $120$ are of the form
$n(15n+3s-5)/2$ or $(3n-s/2)(5n-3+s/2)/2$ with $k\in\mathbb{Z}$. The first identity is proved.

In a similar way, multiplying both sides of \eqref{eq4.3} by \eqref{eq4.2}, we obtain 
\begin{align}
& \bigg(1 + 2 \sum_{n=1}^{k} (-1)^n q^{2n^2} \bigg) \sum_{n=0}^{\infty} \frac{q^{n(n+1)}}{(q;q)_{2n+(s-1)/2}} -\sum_{n=-\infty}^\infty (-1)^{n(n+s)/2}q^{n(5n-s)/2} \nonumber \\
& =2(-1)^k \frac{(-q^2;q^2)_k}{(q^2;q^2)_k} \sum_{n=k+1}^{\infty} \frac{q^{2n(k+1)}(-q^{2(n+1)};q^2)_{\infty}}{(1-q^{2n})(q^{2(n+1)};q^2)_{\infty}}
 \sum_{n=-\infty}^\infty (-1)^{n(n+s)/2}q^{n(5n-s)/2}.\label{Id2}
\end{align}
For each $s\in\{1,3\}$, it is elementary to see that squares that are congruent to $s^2$ modulo $40$ are of the form
$n(5n-s)/2$ with $n\in\mathbb{Z}$. Thus we deduce the second identity.

\section{Proof of Corollary \ref{C2}}
\label{S3}

For each $s\in\{2,4\}$, the generating function for $R_s(n)$ can be written as follows:
$$\sum_{n=0}^\infty R_s(n) q^n 
= \frac{(q^{6-s},q^{4+s},q^{10};q^{10})_\infty}{(q;q)_\infty}.$$
By the Jacobi triple product identity \eqref{eq:JTP}, with $q$ replaced by $q^{10}$ and $z$ replaced by $q^{6-s}$, we get
$$
\sum_{n=-\infty}^\infty (-1)^n q^{n(5n+1-s)} = (q^{6-s},q^{4+s},q^{10};q^{10})_\infty.
$$
Thus we obtain
\begin{align*}
\sum_{n=0}^\infty R_s(n) q^n 
& = \left( \sum_{n=0}^\infty p(n) q^n \right) \left( \sum_{n=-\infty}^\infty (-1)^n q^{n(5n+1-s)} \right) \\
& = \sum_{n=0}^\infty \left( \sum_{k=-\infty}^\infty (-1)^k p\big(n-k(5k+1-s)\big)  \right) q^n.
\end{align*}
Equating the coefficient of $q^n$ in this equation, we obtain the following decomposition of $R_s(n)$ in terms of the partition function $p(n)$:
$$
R_s(n) = \sum_{k=-\infty}^\infty (-1)^k p\big(n-k(5k+1-s)\big).
$$
It is an easy exercise to show that the squares congruent to $(s-1)^2$ modulo $20$ are of the form $k(5k+1-s)$
with $k\in\mathbb{Z}$. 

For each $s\in\{1,3\}$, considering Watson's quintuple product identity \eqref{WQPI}, with $q$ replaced by $q^{10}$ and $z$ replaced by $q^s$, we get
\begin{align*}
\sum_{n=0}^\infty R^*_s(n) q^n 
& = \frac{1}{(q;q)_\infty} \sum_{n=-\infty}^\infty q^{15n^2-(5-3s)n}(1-q^{10n+s})\\
& = \left( \sum_{n=0}^\infty p(n) q^n\right)  \left( \sum_{n=-\infty}^\infty \big(q^{15n^2-(5-3s)n} - q^{15n^2+(5+3s)n+s} \big)\right) .
\end{align*}
Equating the coefficient of $q^n$ in this equation, we obtain
$$R^*_s(n) =\sum_{k=-\infty}^\infty \Bigg(p\Big(n-\big(15k^2-(5-3s)k\big)\Big) - p\Big(n-\big(15k^2+(5+3s)k+s\big)\Big)\Bigg).$$
It is an easy exercise to show that the squares congruent to $(s+1)^2/4$ modulo $15$ are of the form $15k^2-(5-3s)k$ or $15k^2+(5+3s)k+s$ with $k\in\mathbb{Z}$. 

\section{Open problems}
\label{S4}

Related to relations \eqref{Id1} and \eqref{Id2}, we remark that there is a substantial amount of numerical evidence to state the following conjecture.

\begin{conjecture}
	Let $k$ be a positive integer.
	\begin{enumerate}
		\item For $s\in\{2,4\}$,
\begin{align*}
& \bigg(1 + 2 \sum_{n=1}^{k} (-1)^n q^{2n^2} \bigg) \sum_{n=0}^{\infty} \frac{(-q;q)_n q^{n(3n+s-1)/2}}{(q;q)_{2n+1}} \\ 
& \qquad -\sum_{n=-\infty}^{\infty}(-1)^{n((s-1)n-1)/2}\left( q^{n(15n+3s-5)/2}+q^{(3n-s/2)(5n-3+s/2)/2}\right) 
\end{align*}	
	has nonnegative coefficients if $k$ is even and nonpositive coefficients if $k$ is odd.
	\item For $s\in\{1,3\}$, 
	$$
	\bigg(1 + 2 \sum_{n=1}^{k} (-1)^n q^{2n^2} \bigg) \sum_{n=0}^{\infty} \frac{q^{n(n+1)}}{(q;q)_{2n+(s-1)/2}} -\sum_{n=-\infty}^\infty (-1)^{n(n+s)/2}q^{n(5n-s)/2} 
	$$	
	has nonnegative coefficients if $k$ is even and nonpositive coefficients if $k$ is odd.	
	\end{enumerate}
\end{conjecture}

Assuming this conjecture, we derive the following families of linear inequalities involving $R_s(n)$ and $R^*_s(n)$.

\begin{conjecture}
	Let $k$ be a positive integer. 
	\begin{enumerate}
		\item   For $s\in\{2,4\}$,
		\begin{align*}
		(-1)^k \left(R_s(n)+2\sum_{j=1}^k (-1)^k R_s(n-2j^2) - \rho_s(n) \right) \geqslant 0,
		\end{align*}
		where
\begin{equation*}
\rho_s(n) = 
\begin{cases}
(-1)^{m((s-1)m-1)/2}, &\text{if $n=m(15m+3s-5)/2$ or}\\
& \text{\quad $n=(3m-s/2)(5m-3+s/2)/2$, $m\in\mathbb{Z}$,}\\
0, &\text{otherwise.}
\end{cases} 
\end{equation*}
		\item 	For $s\in\{1,3\}$,
		$$
		(-1)^k\left( R^*_s(n)+2\sum_{j=1}^k R^*_s(n-2j^2)-\rho^*_s(n) \right) \geqslant 0,
		$$	
		where
		\begin{equation*}
		\rho^*_s(n) = 
		\begin{cases}
		(-1)^{m(m+s)/2}, &\text{if $n=m(5m-s)/2$, $m\in\mathbb{Z}$,}\\
		0, &\text{otherwise.}
		\end{cases} 
		\end{equation*}
	\end{enumerate}
\end{conjecture}

It would be very appealing to obtain other $q$-series congruences that involve  
other statistical mechanics partition functions appearing in Baxter's solution of the hard-hexagon model.

\bigskip


\end{document}